\documentclass [11pt]{amsart}
\usepackage[english]{babel}
\usepackage{amsfonts}
\usepackage{amssymb}
\usepackage{amsmath}

\newcounter{remark}
\newcounter{theor}
\newtheorem{thm}{Theorem}[section]
\newtheorem{prop}{Proposition}[section]
\newtheorem{conj}{Conjecture}[section]
\newtheorem{lem}{Lemma}[section]
\newtheorem{defn}{Definition}[section]

\newcommand{\de}{\partial}

\newcommand{\C}{\mathbb{C}}
\newcommand{\R}{\mathbb{R}}
\newcommand{\Pp}{\mathbb{P}}
\newcommand{\K}{\mathcal{K}}
\newcommand{\Kb}{\overline{\mathcal{K}}}
\newcommand{\Kns}{\mathcal{K}_{NS}}
\newcommand{\Kbns}{\overline{\mathcal{K}}_{NS}}
\newcommand{\Oo}[1]{\mathcal{O}(#1)}
\newcommand{\Op}[1]{\mathcal{O}_{#1}}

\newcommand{\Z}{\mathbb{Z}}
\newcommand{\Q}{\mathbb{Q}}

\newcommand{\ve}{\varepsilon}

\newcommand{\ti}[1]{\tilde{#1}}

\newcommand{\vol}{\mathrm{Vol}}

\newcommand{\db}{\overline{\partial}}
\newcommand{\Ric}{\mathrm{Ric}}

\newcommand{\diam}{\mathrm{diam}}

\numberwithin{equation}{section}

\begin{document}

\title[Degenerations]{Limits of Calabi-Yau metrics when the K\"ahler class degenerates}
\author{Valentino Tosatti}
\begin{abstract} 
We study the behaviour of families of Ricci-flat K\"ahler metrics on a projective Calabi-Yau manifold
when the K\"ahler classes degenerate to the boundary of the ample cone. We prove that if the limit class is big and nef the Ricci-flat metrics converge smoothly on compact sets outside a subvariety to a limit incomplete Ricci-flat metric. The limit can also be understood from algebraic geometry. 
\end{abstract}

\thanks{Part of this work was carried out while the author was visiting the Morningside Center of Mathematics in Beijing; 
the author is supported in part
by a Harvard Merit Fellowship.}

 \address{Department of Mathematics \\ Harvard University \\ Cambridge, MA 02138}

  \email{tosatti@math.harvard.edu}

\maketitle

\section{Introduction}
Einstein metrics, namely metrics with constant Ricci curvature, have been an important subject
of study in the field of differential geometry since the early days. The solution of the Calabi Conjecture
given by Yau \cite{yaupnas} in 1976 provided a very powerful existence theorem for K\"ahler-Einstein metrics with negative or zero Ricci curvature (the negative case was also done independently by Aubin \cite{aubin}). This produced a number of nonhomogeneous examples of Ricci-flat manifolds. These spaces have been named Calabi-Yau manifolds by the physicists in the Eighties, and have been throughly studied in several different areas of mathematics and physics. Prompted by the physical intuition of mirror symmetry, mathematicians have studied the ways in which Calabi-Yau manifolds can degenerate when they are moving in
families. In general both the complex and symplectic (K\"ahler) structure are changing, and 
the behaviour is not well understood. In this paper we will consider the case when the complex structure is fixed, and so we will be looking at a single compact projective Calabi-Yau manifold. The K\"ahler class is then allowed to vary inside the ample cone. As long as the class stays inside the cone, the corresponding Ricci-flat metrics vary smoothly, but they will degenerate when the class approaches
the boundary of the cone. We will try to understand this degeneration process and see what the limiting space looks like.

To introduce our results, let us fix some notation first.
Let $X$ be a compact projective Calabi-Yau manifold, of complex dimension $n$.
This is by definition a projective manifold such that
$c_1(X)=0$ in $H^2(X,\R)$.  The real N\'eron-Severi space is by definition
$$N^1(X)_\R=(H^2(X,\Z)_{free}\cap H^{1,1}(X))\otimes\R=N^1(X)_\Z\otimes\R,$$
and we assume that $$\dim N^1(X)_\R= \rho(X)>1.$$ This cohomology space
contains $\Kns$ the ample cone, which is open. Its closure $\Kbns$ is the nef cone. Fix a nonzero class
$\alpha\in\Kbns\backslash\Kns$, which exists precisely when $\rho(X)>1$, and a smooth path $\alpha_t:[0,1]\to\Kbns$ such that 
$\alpha_t\in\Kns$ for $t<1$ and $\alpha_1=\alpha$. For any $t<1$ Yau's Theorem \cite{yau1} gives us 
a unique Ricci-flat K\"ahler metric $\omega_t\in\alpha_t$. Fixing a smooth path of reference
metrics in $\alpha_t$, it can be verified that the Ricci-flat metrics $\omega_t$ vary smoothly, as long as $t<1$.
We have the following very natural\\

\noindent{\bf Question 1:}\ \emph{What is the behaviour of the metrics $\omega_t$ as $t\to 1$?}\\

This question has a long history: it is a special case of a problem by Yau \cite{yau2}, \cite{yau3},
where the complex structure is also allowed to vary;  it has been stated explicitly in this form by
McMullen \cite{ctm} and Wilson \cite{wilson}. Physicists have also looked at this question,
roughly predicting the behaviour that we will describe in Theorem \ref{main1} (see e.g. \cite{hw}).
One of the reasons that makes this question interesting is that the Ricci-flat metrics
are not known explicitly, except in very few cases.

A nef class $\alpha\in N^1(X)_\mathbb{R}$ is called {\em big} if $\alpha^n>0.$
Classical results of Anderson \cite{and}, Bando-Kasue-Nakajima \cite{bkn}, Tian \cite{tian} and 
 more recent results of Cheeger-Colding-Tian \cite{cct} give a partial answer to this question when $\alpha$
is big (we will explain this in section \ref{comp}). Our main theorem, which does not rely on the previous results just quoted, gives a satisfactory answer to Question 1 
in this case (see section \ref{algebra} for definitions).

\begin{thm}\label{main1}
Let $X$ be a compact projective Calabi-Yau manifold, and let $\alpha\in N^1(X)_\R$ be a big and nef class that is not ample.
Then there exist a proper analytic subvariety $E\subset X$ and a smooth incomplete Ricci-flat 
K\"ahler metric $\omega_1$ on $X\backslash E$, that depend only 
on $\alpha$, such that
 for any smooth path $\alpha_t\in\Kns$ with $\alpha_1=\alpha$, the Ricci-flat metrics $\omega_t\in\alpha_t$ converge to $\omega_1$ in the $C^\infty$ topology on compact
sets of $X\backslash E$. 
Moreover $\omega_1$ extends to a closed positive current with continuous potentials
on the whole of $X$, which lies in $\alpha$, and which is the pullback of a singular
Ricci-flat K\"ahler metric on a Calabi-Yau model of $X$ obtained from the contraction map of $\alpha$. 
If $\alpha\in N^1(X)_\Z$, that is if
$\alpha=c_1(L)$ for some line bundle $L$, then $E$ is the null locus of $L$.
\end{thm}

There are many interesting concrete examples of our theorem, and we will examine a few of them in section \ref{examples}.
Roughly speaking, the case when $\alpha$ is nef and big corresponds to a ``non-collapsing'' sequence of metrics,
meaning that the Gromov-Hausdorff limit has the same dimension. The ``collapsing'' case, when $\alpha$ is nef but not 
big, is much harder and we will briefly discuss it at the end of the paper.
We state and prove our results for a path of classes $\alpha_t$, but it's immediate to see that the same
results hold if instead we look at a sequence of classes $\alpha_i$ that converge to $\alpha$.
On the other hand, our result doesn't say anything about the case when the classes $\alpha_t$
approach the boundary of the ample cone without converging to a limiting class, but moving out
to infinity in $N^1(X)_\R$. This case is relevant for mirror symmetry, as it should sometimes be the mirror of a large complex structure limit. Finally let us remark that the projectivity assumptions are only technical, and that we expect that a similar result holds when $X$ is just assumed to be 
K\"ahler, and the ample cone is replaced by the K\"ahler cone (see section \ref{direc}).\\

The paper is organized as follows. In section 2 we recall some definition and results from algebraic geometry. In section 3 we prove a uniform diameter bound and we
compare our results with previous literature. In section 4 we prove our main Theorem \ref{main1}. We use a new Moser iteration argument to get uniform bounds, using the diameter bound from section 3. In section 5 we give some examples where our results apply, and recover in particular a result of Kobayashi and Todorov \cite{kt}. Finally in section 6 we discuss some further directions for research.\\

{\bf Acknowledgments.}\ I would like to thank my advisor Prof. Shing-Tung Yau for suggesting this problem and for constant support. I also thank Prof. Curt McMullen and Prof. M.S. Narasimhan for inspiring conversations, Chen-Yu Chi, Jian Song, G\'abor Sz\'ekelyhidi and Ben Weinkove for useful comments, and the referee for some interesting suggestions.

\section{Some facts from algebraic geometry}\label{algebra}
In this section we will review some definitions and results from algebraic geometry, mainly from Mori's Program, that will be used in the proof.\\

Let $X$ be a compact Calabi-Yau $n$-fold, that is a compact K\"ahler manifold of dimension $n$ and such that 
$c_1(X)=0$ in $H^2(X,\R)$. We don't insist that $X$ be simply connected. 
Notice that it follows that $aK_X\cong \Op{X}$ for some integer $a>0$: in fact by Theorem 1 in \cite{beauville}
a finite unramified $a:1$ cover of $X$, $p:\ti{X}\to X$, has trivial canonical bundle. But we have that
$p^*K_X\cong K_{\ti{X}}\cong\Op{\ti{X}}$ and then Lemma 16.2 in \cite{bpv} implies that $aK_X\cong \Op{X}$.
This can be rewritten as $K_X\sim_\Q 0$ where $\sim_\Q$ indicates $\Q$-linear equivalence of Cartier $\Q$-divisors.
For the rest of this section we will assume that $X$ is projective.

\begin{defn} A projective variety $X$ has canonical singularities if it is normal, if $rK_X$ is Cartier 
for some $r\geq 1$ and if there exists a resolution $f:Y\to X$ such that
$$rK_Y=f^*(rK_X)+\sum_i a_i E_i,$$
where $E_i$ ranges over all exceptional prime divisors of $f$, and $a_i\geq 0$.
\end{defn}

\begin{defn}[Wilson \cite{wilson1}]A \emph{Calabi-Yau model} $Y$ is a normal projective variety with canonical singularities
and such that $K_Y\sim_\Q 0$.
\end{defn}

Let $L$ be a nef line bundle on $X$, and let $\kappa(X,L)$ be its Iitaka dimension, that is
$$\kappa(X,L)=m\quad \iff \quad h^0(X,kL)\sim k^m \textrm{ for all } k \textrm{ large enough}$$
and $\kappa(X,L)=-\infty$ if $kL$ has no sections for all $k\geq 0$. We call $\nu(X,L)$ its numerical dimension, that is
the largest nonnegative integer $m$ such that there exists an $m-$cycle $V$ such that $(L^m\cdot V)>0$.
It is always true that $$\kappa(X,L)\leq\nu(X,L)\leq n.$$
\begin{defn} If $\kappa(X,L)=\nu(X,L)$ we say that $L$ is \emph{good} (or abundant). If the complete linear system $|kL|$ is base-point-free for some $k\geq 1$ 
we say $L$ is \emph{semiample}.
\end{defn}

When $|kL|$ is base-point-free, we get a morphism $\Phi_{|kL|}:X\to \mathbb{P}H^0(X,kL)^*$
that satisfies $kL=\Phi_{|kL|}^*\Oo{1}$.
Notice that if $L$ is \emph{big}, that is $\kappa(X,L)=n$, then it is automatically good. 
The following is an immediate consequence of the base-point-free Theorem (Theorem 6.1.11 in \cite{kmm}).

\begin{thm}[Kawamata]\label{bpf} Assume $X$ is a projective Calabi-Yau. If $L$ is good then it is semiample.
\end{thm}

The next theorem is classical (see Theorem 2.1.27 in \cite{laza}).
\begin{thm}[Iitaka]\label{semiample}Let $L$ be semiample. Then there exists
a surjective morphism $f:X\to Y$ where $Y$ is a normal irreducible variety, $f_*\Op{X}=\Op{Y}$, and $L=f^*A$ for
some ample line bundle $A$ on $Y$. In fact $f=\Phi_{|kL|}$ for all $k$ sufficiently divisible.
\end{thm}
We will call $f$ the {\em contraction map} of $L$. There is a version of the base-point-free Theorem for Cartier $\R$-divisors, essentially due to Shokurov \cite{shok}.
If $D$ is a Cartier $\R$-divisor on $X$ we say that $D$ is semiample if there exist $Y$ a normal irreducible projective variety, $f:X\to Y$ a surjective morphism with $f_*\Op{X}=\Op{Y}$, and $A$ an ample $\R$-divisor on $Y$ such that $D\sim_{\R} f^*A$. Again we will call $f$ the contraction map of $D$. Then the following holds (Theorem 7.1 in \cite{hm}): 

\begin{thm}\label{rsemi} Assume $X$ is a projective Calabi-Yau. If $D$ is a Cartier $\R$-divisor which is nef and big, 
then it is semiample.
\end{thm}

The contraction map of $D$ is in fact also the contraction map of a suitable nef and big line bundle $L$
(see the proof of Proposition \ref{semipositive}). We also have the following 
theorem (Theorem 5.7 in \cite{kawa} or Theorem 1.9 in \cite{kawa2}).

\begin{thm}[Kawamata]\label{poly} Assume $X$ is a projective Calabi-Yau. Then the subcone of $\Kbns$ given by nef and big classes
is locally rational polyhedral.
\end{thm}

If $L$ is a line bundle, its stable base locus is the intersection of the base loci of $|mL|$ for all $m\geq 1$. It is equal to the base locus of $|mL|$ for some $m$ (see Prop. 2.1.21 in \cite{laza}). If $L$ now is nef and big, we define the augmented base locus of $L$, $\mathbf{B}_+(L)$, to be the stable base locus of $L-\ve H$ for any $H$ ample divisor and any $\ve>0$ small enough rational number.
This definition is well-posed (see Lemma 10.3.1 in \cite{laza}) and a theorem of Nakamaye (\cite{naka}, \cite{laza}) says that $\mathbf{B}_+(L)$ is equal to the null locus of $L$, that is the union of all
positive-dimensional subvarieties $V\subset X$ such that $(L^{\dim V}\cdot V)=0$.\\

Finally let us state a well-known conjecture (see 10.3 of Peternell's lectures in \cite{mp}).
\begin{conj}\label{bpf2} Assume $X$ is a projective Calabi-Yau. If $L$ is a nef line bundle, then $L$ is semiample.
\end{conj}
If $L$ is effective, this conjecture follows from the log abundance conjecture. Indeed for any small rational $\ve>0$, the
pair $(X,\ve L)$ is klt, and the log abundance conjecture would imply that $K_X+\ve L\sim_{\Q}\ve L$ is semiample.

Notice that when $X$ is a surface, Conjecture \ref{bpf2} holds: in fact if $L$ is nef and non trivial,
then $H^2(X,L)=H^0(X,K_X-L)=0$ and by Riemann-Roch
$$\dim H^0(X,L)\geq 2+\frac{1}{2}L\cdot L\geq 2,$$
thus $L$ is effective. Then we can apply the log abundance theorem for surfaces (see e.g. \cite{log}) and get the result. 

\section{Preliminary remarks}\label{comp}
In this section we will prove a uniform diameter bound and use this to compare our results to the existing literature. The diameter bound is valid in general, without any projectivity or bigness assumptions.\\

Let the setting be as in the Introduction, namely let $X$ be a projective Calabi-Yau $n$-fold
and $\alpha\in N^1(X)_\R$ a big and nef class that is not ample. Given $\alpha_t:[0,1]\to\Kbns$ a smooth path such that $\alpha_t\in\Kns$ for $t<1$ and $\alpha_1=\alpha$, Yau's Theorem \cite{yau1} gives us 
a unique Ricci-flat K\"ahler metric $\omega_t\in\alpha_t$ for any $t<1$.
Then $$\sqrt{-1}\de\db\log\frac{\omega_t^n}{\omega_0^n}=\Ric(\omega_0)-\Ric(\omega_t)=0$$
implies that $\omega_t^n=B_t\omega_0^n$ for some constant $B_t$, which is easily computed by
\begin{equation}\label{dame3}
\alpha_t^n=\int_X\omega_t^n=B_t\int_X \omega_0^n=B_t\alpha_0^n.
\end{equation}
In particular $B_t>0$ and $0<\lim_{t\to 1}B_t<\infty$, which means that the volume form of
$\omega_t$ is uniformly equivalent to the volume form of $\omega_0$. The main theorem of this section is the following

\begin{thm}\label{diameter}
Let $(X,\omega_0)$ be a compact $n$-dimensional Ricci-flat K\"ahler manifold and let $\omega$ be another Ricci-flat K\"ahler metric such that
\begin{equation}\label{dame2}
\int_X \omega_0^{n-1}\wedge\omega\leq c_1,
\end{equation}
for some constant $c_1$. Then the diameter of $(X,\omega)$ is
bounded above by a constant that depends only on $n,c_1,\omega_0$.
\end{thm}
Applying this to $\omega=\omega_t$ we see that the diameter of $(X,\omega_t)$
is uniformly bounded as $t$ approaches $1$.
To prove Theorem \ref{diameter} we need a lemma, which appears as Lemma 1.3 in \cite{dps}. For the reader's convenience, we
include a proof here.
\begin{lem}\label{dps1}
In the above situation there exists a constant $C_1$ that depends only on $n,c_1,\omega_0$, such that given any $\delta>0$ there exists 
an open set $U_{\delta}\subset X$ such that its diameter with respect to $\omega$ is less
than $C_1\delta^{-1/2}$ and its volume with respect to $\omega_0$ is at least $\int_X\omega_0^n-\delta$. 
\end{lem}
\begin{proof}
First notice that 
\eqref{dame2}
gives a uniform $L^1$ bound on $\omega$.
Up to covering $X$ by finitely many charts, we may assume that $X=K$ is a compact convex set in $\C^n$, and we will denote by $g_E$ the Euclidean metric on $K$. If $x_1, x_2\in K$, we denote by $[x_1,x_2]$ the segment joining them in $K$, and we compute the average of the length square of $[x_1,x_2]$ with respect to $\omega$, when the endpoints vary. Using Fubini's Theorem and the Cauchy-Schwarz inequality we get
\begin{equation}\label{length1}
\begin{split}
&\int_{K\times K}\left(\int_0^1 \sqrt{\omega((1-s)x_1+sx_2)(x_2-x_1)}ds\right)^2 dx_1 dx_2\\
&\leq \|x_2-x_1\|^2_{g_E}\int_0^1 \int_{K\times K}|\omega((1-s)x_1+sx_2)|dx_1dx_2 ds\\
&\leq \diam_{g_E}^2(K) 2^{2n}\biggl(\int_0^{\frac{1}{2}} \int_{K\times K}|\omega(y+sx_2)|dydx_2 ds\\
&+ \int_{\frac{1}{2}}^1 \int_{K\times K}|\omega((1-s)x_1+y)|dydx_1 ds\biggr)\\
&\leq  \diam_{g_E}^2(K) 2^{2n}\vol_{g_E}(K)\|\omega\|_{L^1(K)}\leq C_1,
\end{split}
\end{equation}
where $C_1$ is a uniform constant, we changed variable $y=(1-s)x_1$ if $s\leq \frac{1}{2}$ and $y=sx_2$ when $s\geq \frac{1}{2}$ and
integrated first with respect to $y$. Then the set $S$ of pairs $(x_1,x_2)\in K\times K$ such that the length of $[x_1,x_2]$ with respect to $\omega$ is more than $(C_1/\delta)^{1/2}$ has Euclidean measure less than or equal $\delta$: otherwise 
\begin{equation*}\label{length2}
\begin{split}
&\int_{K\times K}\left(\int_0^1 \sqrt{\omega((1-s)x_1+sx_2)(x_2-x_1)}ds\right)^2 dx_1 dx_2\\
&\geq \int_{S}\left(\int_0^1 \sqrt{\omega((1-s)x_1+sx_2)(x_2-x_1)}ds\right)^2 dx_1 dx_2\geq
\frac{C_1}{\delta}\vol_{g_E}(S)
\end{split}
\end{equation*}
which is more than $C_1$, and this contradicts \eqref{length1}. If $x_1\in K$ we let $S(x_1)$ to be the set of the $x_2\in K$ such that $(x_1,x_2)\in S$, and we let $Q$ to be the set of the $x_1\in K$ such that
$\vol_{g_E}(S(x_1))\geq\frac{1}{2}\vol_{g_E}(K)$ and $R$ to be the set of $(x_1,x_2)\in S$ such that
$x_1\in Q$. Then by Fubini's Theorem
$$\delta\geq \vol_{g_E}(R)=\int_R dx_2 dx_1=\int_Q \left(\int_{S(x_1)}dx_2\right) dx_1\geq
\frac{1}{2}\vol_{g_E}(K)\vol_{g_E}(Q),$$
and so $\vol_{g_E}(Q)\leq \frac{2\delta}{\vol_{g_E}(K)}.$ We let $U_{\delta}=K\backslash Q$.
Then $U_{\delta}$ is open and if $x_1,x_2\in U_{\delta}$ then $\vol_{g_E}(S(x_i))<\frac{1}{2}\vol_{g_E}(K)$, for $i=1,2$. Hence $\vol_{g_E}((K\backslash S(x_1))\cap (K\backslash S(x_2)))>0$ and so this set is nonempty. If $y$ belongs to it, then $(x_1,y)$ and $(x_2,y)$ are not in $S$, which means that the lengths with respect to $\omega$ of the segments $[x_1,y]$ and $[y,x_2]$ are both less than 
$(C_1/\delta)^{1/2}$. Concatenating these two segments we get a path from $x_1$ to $x_2$ with length less than $2(C_1/\delta)^{1/2}$. We also have that
$$\vol_{\omega_0}(Q)\leq C_2 \vol_{g_E}(Q)\leq \frac{2C_2\delta}{\vol_{g_E}(K)}.$$
Up to adjusting the constants, this is what we want.
\end{proof}

\begin{proof}[Proof of Theorem \ref{diameter}]

Choose $\delta\leq\min(C_1^2,1/2\int_X\omega_0^n)$,
and pick any $p\in U_{\delta}$. If we denote the metric ball of $\omega$ centered at
$p$ and with radius $r$ by $B(p,r)$, then we get that $U_{\delta}\subset B(p,C_2)$, where
$C_2=C_1\delta^{-1/2}\geq 1$. Then $$\int_{B(p,C_2)}\omega_0^n
\geq \int_{U_{\delta}}
\omega_0^n\geq\frac{1}{2}\int_X\omega_0^n.$$
Proceeding as in \eqref{dame3} we see that $\omega^n=B\omega_0^n$, where 
\begin{equation}\label{dami}
B=\frac{\int_X\omega^n}{\int_X\omega_0^n}. 
\end{equation}
So we get that 
\begin{equation}\label{volume}
\int_{B(p,C_2)}\omega^n\geq B C_3,
\end{equation}
for some constant $C_3>0$ independent of $\omega$. Since $\Ric(\omega)=0$, the Bishop volume comparison Theorem
and \eqref{volume} give that 
\begin{equation}\label{volume2}
\int_{B(p,1)}\omega^n\geq \frac{\int_{B(p,C_2)}\omega^n}{C_2^{2n}}\geq B C_4>0.
\end{equation}
The following lemma is due to Yau (see e.g. Theorem I.4.1 in \cite{sy}), but we provide a proof for completeness.
\begin{lem}\label{sy1}
Let $(M^{2n},g)$ be a closed Riemannian manifold with $\Ric(g)\geq 0$, let $p\in M$ and $1<R<\diam(X,g)$. Then
$$\frac{R-1}{4n}\leq\frac{\vol(B(p,2(R+1)))}{\vol(B(p,1))}.$$
\end{lem}
\begin{proof}
Choose $x_0\in\partial B(p,R)$, so that $d(x_0,p)=R$, and denote by $\rho(x)=d(x,x_0)$.
The Laplacian comparison theorem gives $\Delta\rho^2\leq 4n$ in the sense of distributions.
Let $\varphi(x)=\psi(\rho(x))$ where 
$$\psi(t)= \left\{ \begin{array}{lll}
         1 & \mbox{if $0 \leq t\leq R-1$},\\
        \frac{1}{2}(R+1-t) & \mbox{if $R-1<t<R+1$},\\
         0 & \mbox{if $t\geq R+1$}.\end{array} \right.$$
Then $\varphi$ is a nonnegative Lipschitz function supported in $B(x_0,R+1)$, and we have that
\begin{equation*}\begin{split}
\int_M\varphi\Delta\rho^2dV_g&=-\int_{B(x_0,R+1)}\nabla\varphi\cdot\nabla\rho^2dV_g=-2\int_{B(x_0,R+1)}\rho|\nabla\rho|^2\psi'(\rho(x))dV_g\\
&=\int_{B(x_0,R+1)\backslash B(x_0,R-1)}\rho dV_g\\
&\geq (R-1)\vol(B(x_0,R+1)\backslash B(x_0,R-1)),
\end{split}\end{equation*}
and also
$$\int_M\varphi\Delta\rho^2dV_g\leq 4n\int_{B(x_0,R+1)}\varphi dV_g\leq 4n\vol(B(x_0,R+1)).$$
Notice that $B(p,1)\subset B(x_0,R+1)\backslash B(x_0,R-1)$ and so the previous two equations give
$$(R-1)\vol(B(p,1))\leq 4n\vol(B(x_0,R+1)).$$
The conclusion follows from the fact that $B(x_0,R+1)\subset B(p,2(R+1))$.
\end{proof}
Lemma \ref{sy1} gives that for any $1<R<\diam(X,\omega)$ we have
\begin{equation}\label{volume3}
\frac{R-1}{4n}\leq\frac{\int_{B(p,2(R+1))}\omega^n}{\int_{B(p,1)}\omega^n}.
\end{equation}
Choosing $R=\diam(X,\omega)-1$ and using \eqref{volume2}, \eqref{dami} we get
$$\diam(X,\omega)\leq 2+\frac{4n}{B C_4}\int_X\omega^n= 2+\frac{4n}{C_4}\int_X\omega_0^n,$$
which is bounded independent of $\omega$. This completes the proof of Theorem \ref{diameter} (a somewhat similar argument can be found in \cite{paun1}). 
\end{proof}

Once we have the diameter bound, we can apply the Bishop volume comparison Theorem again and get
that for any point $p\in X$ and any $r>0$, $t<1$
\begin{equation}\label{volume4}
\int_{B_t(p,r)}\omega_t^n\geq r^{2n} \frac{\int_X\omega_t^n}{\diam(X,\omega_t)^{2n}}\geq c r^{2n},
\end{equation}
where $c>0$ is a uniform constant.
A well-known computation in Chern-Weil theory gives
$$\frac{1}{n(n-1)}\int_X \|\textrm{Rm}_t\|^2_t\omega_t^n=\int_X c_2(X,\omega_t)\wedge\omega_t^{n-2}=c_2(X)\cdot\alpha_t^{n-2}\leq C,$$
where $\textrm{Rm}_t$ is the Riemann curvature tensor of $\omega_t$ and $c_2(X,\omega_t)$ is the
second Chern form of $\omega_t$. If $n=2$ we can thus apply Theorem C of \cite{and}, Theorem 5.5 of \cite{bkn}
or Proposition 3.2 of \cite{tian}
and get that a subsequence of $(X,\omega_t)$ converges to an  Einstein orbifold
with isolated singularities in the Gromov-Hausdorff topology, and also in the $C^\infty$ topology on compact sets outside the orbifold points. If $n>2$ these theorems require a uniform
bound on $$\int_X \|\textrm{Rm}_t\|^n_t\omega_t^n,$$
which in general can not be expressed in terms of topological data as above. 
Instead when $n>2$ we apply a general theorem of Gromov \cite{gromov} that says that any sequence of compact Riemannian manifolds of dimension $2n$
with diameter bounded above and Ricci curvature bounded below, has a subsequence that converges in
the Gromov-Hausdorff topology to a 
compact length space. Thus a subsequence of $(X,\omega_t)$ converges to a compact metric space $Y$, and Theorem 1.15 in \cite{cct} says that $Y$ is a complex manifold
outside a rectifiable set $R\subset Y$ of real Hausdorff codimension at least $4$. Moreover their Theorem 9.1
gives supporting evidence that $R$ should in fact be a complex subvariety of $Y$.

On the other hand our Theorem \ref{main1} gives the convergence of the whole sequence of metrics, and not just of a subsequence, and the limit metric is uniquely determined by the class $\alpha$.
When $n>2$ the convergence we get is stronger than Gromov-Hausdorff convergence, but it only happens outside the singular set $E$. Also we see precisely who
the limit space $Y$ is, namely the Calabi-Yau model of $X$ obtained from the contraction map of $\alpha$. It has canonical singularities,
so its singular set is a subvariety
of complex codimension at least $2$, and when $n=2$ canonical singularities
are precisely rational double points, that are of orbifold type. We will discuss the case $n=2$ with more
details in section \ref{examples}.
Let us also mention the results of Ruan \cite{ruan}. He studies the Gromov-Hausdorff limits of sequences
of K\"ahler metrics on a fixed compact manifold $X$, with uniformly bounded sectional curvature. Roughly speaking, he proves that there exists
an analytic subvariety $E\subset X$ such that a subsequence of the metrics converges in the Gromov-Hausdorff topology on $X\backslash E$ to $\omega$ a smooth Hermitian form, which is either K\"ahler (non-collapsing)
or pointwise nonnegative with determinant zero (collapsing). Moreover in the collapsing case, the kernel of
$\omega$ gives a holomorphic foliation with singularities on $X$. Unfortunately in our setting the curvature is not bounded in general, so Ruan's results don't apply, but our Theorem \ref{main1} gives in particular Ruan's conclusion in the non-collapsing case. We will discuss the collapsing case in section \ref{direc}.

Of course, the above-mentioned results apply in more general situations than ours. Also, all the results in this section work in the case when $X$ is not projective, and the ample cone is replaced by the (bigger) K\"ahler cone. Then the above theorems still apply, but for technical reasons our Theorem \ref{main1} doesn't (see section \ref{direc} for more discussions).

\section{Limits of Ricci-flat metrics}
In this section we will prove Theorem \ref{main1}. The idea is to carefully set up a family 
of complex Monge-Amp\`ere equations that degenerate in the limit, and prove estimates for the solutions
that are uniform outside a subvariety.\\

We begin with a

\begin{prop}\label{semipositive}
Let $X$ be a projective Calabi-Yau $n$-fold, and $\alpha\in N^1(X)_\R$ a big and nef class that is not ample.
Then there exists $\omega\in\alpha$ a smooth real $(1,1)$ form that is {\em pointwise} nonnegative and which is K\"ahler outside a proper analytic subvariety of $X$.
Moreover if $\alpha_t:[0,1]\to\Kbns$ is a smooth path such that $\alpha_t\in\Kns$ for $t<1$ and $\alpha_1=\alpha$,
then we can find a continuous family of K\"ahler forms $\beta_t\in\alpha_t$, $t<1$, such that $\beta_t\to\omega$
in the $C^\infty$ topology as $t$ approaches $1$.
\end{prop}
\begin{proof} 
Let's assume first that that $\alpha=c_1(L)$ for some line bundle $L$, which is equivalent to
requiring that $\alpha\in N^1(X)_\Z$.
Now $L$ is nef and big and so Theorem \ref{bpf} implies that $L$ is semiample, so there exists some $k\geq 1$ 
such that $kL$ is globally generated. This gives a morphism $f:X\to\Pp^N$ such that $f^*\Oo{1}=kL$. If we
let $\omega_{FS}$ be the Fubini-Study metric on $\Pp^N$, then 
$\omega=\frac{f^*\omega_{FS}}{k}$ is a pointwise nonnegative smooth real $(1,1)$ form
in the class $\alpha$. Moreover $\omega$ is K\"ahler outside the exceptional set of $f$, which is a
proper subvariety of $X$.
If $\alpha\in N^1(X)_\Q$, then $k\alpha\in N^1(X)_\Z$
for some integer $k\geq 1$, and we can proceed as above.
If finally $\alpha\in N^1(X)_\R$ then by Theorem \ref{poly} we know
that the subcone of nef and big classes is locally rational polyhedral.
Hence $\alpha$ lies on a face of this cone which is cut out by linear equations with rational coefficients.
It follows that rational points on this face are dense, and it is then
possible to write $\alpha$ as a linear combination of classes in $N^1(X)_\Q$ which are nef and big, 
with nonnegative coefficients.
It is now clear that we can represent $\alpha$ by a smooth nonnegative form $\omega$.
Notice that all of these classes give the same contraction map $f:X\to Y$, because they lie on the same
face. This map is then also the contraction map of $\alpha$, and $\omega$ is again K\"ahler outside
the exceptional set of $f$. 

Now fix a ball $\mathcal{U}$ in 
$N^1(X)_\R$ centered at $\alpha$, such that $\Kns\cap\mathcal{U}$ is defined by
$\{\Phi_\beta>0\}_{1\leq\beta\leq k}$ where the $\Phi_\beta$ are linear forms with rational coefficients. Since the big cone is open, up to shrinking $\mathcal{U}$ we may also assume that all the classes in $\partial\Kns\cap\mathcal{U}$ are big. We may add some more linear forms to the $\Phi_\beta$,
until they define a strongly convex rational polyhedral cone $C$ which is contained in $\Kbns\cap\mathcal{U}$. We can then write
$$C=\left\{\sum_{i=1}^\ell a_i\gamma_i\ \bigg|\ a_i\geq 0\right\},$$
where the $\gamma_i$ are nef and big classes in $\mathcal{U}$.
We claim that, when $t$ is bigger than some $t_0<1$, it is possible to write the path $\alpha_t$ as $\sum_i a_i(t)\gamma_i$
where the functions $a_i(t)$ are continuous and nonnegative.
Assume first that the cone $C$ is simplicial, which means that the $\gamma_i$ are linearly independent. Then the path $\alpha_t$ enters and eventually stays in $C$, and so it can be expressed uniquely as
\begin{equation}\label{scrittura}
\alpha_t=\sum_{i=1}^\ell a_i(t)\gamma_i,
\end{equation}
where the $a_i(t)$ are smooth and nonnegative, $t_0\leq t\leq1$. If on the other hand $C$ is not simplicial,
it can be written as a finite union of simplicial subcones that intersect only along faces, and that
are spanned by some linearly independent
subsets of the $\gamma_i$. On any time interval when $\alpha_t$ belongs to the interior of
a simplicial cone, the coefficients $a_i(t)$ in \eqref{scrittura} vary smoothly, and on a common face of 
two simplicial cones the coefficients agree, hence the $a_i(t)$ vary continuously when $t_0\leq t<1$. Moreover since we only have finitely many simplicial subcones, we see that as $t\to 1$ the $a_i(t)$ converge to
the coefficients of $\alpha_1$ in any of the simplicial cones that contain it, and so the $a_i(t)$ are continuous on the whole interval $t_0\leq t\leq 1$.

By the first part of the proof we know that we can choose $\delta_i\in\gamma_i$ a smooth nonnegative
representative, for all $i$.
Choose a smooth function $\ve(t):[t_0,1]\to\R$ that is positive on $[t_0,1)$ and $\ve(1)=0$, and that is
small enough so that the classes $\ti{\alpha}_t=\alpha_t-\ve(t)\alpha_{t_0}$ are ample for
all $t_0\leq t<1$. Then the new path $\ti{\alpha}_t$
is also converging to $\alpha$ as $t\to 1$, and by the previous claim we can write
$$\ti{\alpha}_t=\sum_{i=1}^\ell \ti{a}_i(t)\gamma_i,$$
where $\ti{a}_i(t)$ is a continuous nonnegative function, for all $i$.
Then the smooth $(1,1)$ forms
$$\ti{\beta}_t:=\sum_{i=1}^\ell \ti{a}_i(t)\delta_i$$
are nonnegative representatives of $\ti{\alpha}_t$ that vary continuously in $t$. When $t$ approaches $1$, the forms
$\ti{\beta}_t$ converge in the $C^\infty$ topology to a smooth nonnegative form $\ti{\omega}$ representing $\alpha$.
If $\chi$ is a K\"ahler form in $\alpha_{t_0}$, then the forms $\beta_t=\ti{\beta}_t+\ve(t)\chi$
defined on $[t_0,1)$ are K\"ahler, represent $\alpha_t$ and converge to $\ti{\omega}$ as $t\to 1$.
Up to replacing $\omega$ by $\ti{\omega}$, this gives the desired family of forms on $[t_0,1)$.
It is very easy to extend the family $\beta_t$ on the whole $[0,1)$,
and since we're not going to use this, we leave the proof to the reader.
\end{proof}

Of course, a similar statement holds if we are given a sequence of ample classes $\alpha_i$ converging
to $\alpha$, instead of a path.

Let us now recall some notation and facts from analytic geometry.
If $X$ is any complex manifold and $\omega$ is a Hermitian form on $X$, we will denote by
$PSH(X,\omega)$ the set of all upper semicontinuous (usc) functions $\varphi:X\to [-\infty,+\infty)$ 
such that $\omega+\sqrt{-1}\de\db\varphi$ is a positive current. In the case when $(X,\omega)$ is K\"ahler, 
then all K\"ahler potentials for $\omega$ belong to $PSH(X,\omega)$. A fundamental result by
Bedford-Taylor \cite{bt} says that the Monge-Amp\`ere operator $(\omega+\sqrt{-1}\de\db\varphi)^n$ is
well defined whenever $\varphi\in PSH(X,\omega)$ is locally bounded. Let's also recall the definition
of a singular K\"ahler metric \cite{eyss} on a (possibly singular) algebraic variety $X$. This is given by specifying
its K\"ahler potentials on an open cover $(U_i)$ of $X$, that are usc functions $\varphi_i:U_i\to [-\infty,+\infty)$
with the following property: $\varphi_i$ extends to a plurisubharmonic function on an open set $V_i\subset \mathbb{C}^m$
where $U_i\subset V_i$ is a local embedding. We refer the reader to section 7 of \cite{eyss} for the definition of a singular Ricci-flat K\"ahler metric and for a proof that they always
exist on Calabi-Yau models. With these facts in mind, we can now give the

\begin{proof}[Proof of Theorem \ref{main1}]
Proposition \ref{semipositive} gives us $\omega\in\alpha$ a smooth nonnegative representative,
and $\beta_t\in\alpha_t$ continuously varying K\"ahler forms, when $t<1$, such that $\beta_t\to\omega$ as $t\to 1$. 
As in the proof of Proposition \ref{semipositive}, there is a contraction map  $f:X\to Y$ such that 
$Y$ is a normal irreducible projective variety, $f$ is birational and $f_*\Op{X}=\Op{Y}$.
Moreover $\omega$ is the
pullback of a (singular) K\"ahler metric on $Y$, and it is K\"ahler outside the exceptional set of $f$.
Then setting $D_0=0$ as Cartier divisors on $Y$, we have $aK_X=f^*D_0$ for some
integer $a>0$, so
$$f_*(aK_X)=D_0=0$$
holds as Weil divisors, but since $f$ is birational we also have $f_*(aK_X)=aK_Y$ (as Weil divisors),
hence $aK_Y$ is Cartier and is equal to zero. So we have $f^*K_Y=K_X$ as $\Q$-divisors, which implies
that $Y$ has at most canonical singularities and is a Calabi-Yau model
(see also Corollary 1.5 of \cite{kawa}).

Denote by $\Omega$ the smooth volume form on $X$ given by $$\Omega=\frac{\omega_0^n}{\int_X\omega_0^n},$$
which satisfies $\int_X \Omega=1$.
We can write $\Omega=F\omega^n,$
where $F\in L^1(\omega^n)$, $F>0$. The following argument to show that
actually $F\in L^p(\omega^n)$ for some $p>1$ is similar to Lemma 3.2 in \cite{eyss}.
First of all $1/F$ is smooth, nonnegative, and vanishes precisely on the exceptional
set of $f$. Fixing local coordinates $(z^i)$ on a polydisc $D\subset X$ and a local embedding $G:f(D)
\to \mathbb{C}^m$, we see that $1/F$ is comparable to
$$\left|\frac{\de G}{\de z^1}\wedge\dots\wedge\frac{\de G}{\de z^n}\right|^2$$
on $D$. But this is in turn comparable to $$\sum_{i=1}^r|g_i|^2,$$
where the $g_i$ are holomorphic functions on $D$, and so $F^\ve\in L^1(D,\Omega)$
for some small $\ve>0$ that depends on the vanishing orders of the $g_i$. Then
\begin{equation}\label{integral}
\int_D F^{1+\ve}\omega^n=\int_D F^\ve \Omega<\infty.
\end{equation} 
The compactness of $X$ gives $F\in L^{1+\ve}(\omega^n)$, and so
we can apply Theorem 2.1 and Proposition 3.1 of \cite{eyss} (which rely on the seminal work of Ko\l odziej \cite{kol}) to get
a unique bounded $\varphi\in PSH(X,\omega)$ such that
\begin{equation}\label{eqn1}
(\omega+\sqrt{-1}\de\db\varphi)^n=\alpha^n\Omega,
\end{equation}
and $\sup_X\varphi=0$. We then embed $Y$ into projective space and extend $\omega$
to a K\"ahler form in a neighborhood of $Y$ as in Proposition 3.3 of \cite{dp}.
Composing the embedding with $f$ we get a morphism which is birational with the image,
with connected fibers, and we can then apply Theorem 1.1 in \cite{zhang}
(see also \cite{zhangthesis} and Remark 5.2 in \cite{dz}) and get that $\varphi$ is continuous.
Moreover we can see that $\varphi$ descends to a function
on $Y$: if $V$ is a fiber of $f$, the restriction of $\varphi$ to $V$ is a plurisubharmonic function,
because $\omega|_V=0$. Desingularizing $V$ and applying the maximum principle we see that
$\varphi|_V$ has to be constant, and so $\varphi$ descends to $Y$.
 Since $\omega$ by construction is the pullback of a (singular) K\"ahler form on $Y$, we
see that $\omega+\sqrt{-1}\de\db\varphi$ is a singular Ricci-flat metric on $Y$, in the terminology of \cite{eyss}.
On $X$, the closed positive current $\omega_1=\omega+\sqrt{-1}\de\db\varphi$ clearly lies in the class $\alpha$ and has continuous potentials.

Intuitively, our goal is to get estimates in the open set where $\omega$ is positive. This can be done rigorously in the following way,
which was first used by H.Tsuji \cite{tsuji} (see also \cite{tz}, \cite{lanave} for a recent revisiting of his approach).
Since $\alpha$ is nef and big, by Kodaira's lemma (Example 2.2.23 in \cite{laza}) there exists $E$ an effective Cartier $\R$-divisor such that
for all $\ve>0$ small enough, $\alpha-\ve E=\kappa_\ve$ is K\"ahler. 

We will show that $\varphi$ is smooth on $X\backslash E$, and so $\omega_1$ is a smooth Ricci-flat
metric there, and that the Ricci-flat metrics $\omega_t$ converge to $\omega_1$ in
the $C^\infty$ topology on compact sets of $X\backslash E$. Notice that
the metric $\omega_1$ on $X\backslash E$ cannot be complete, since its diameter is finite by Theorem \ref{diameter}. Once this is proved, 
we can repeat the argument for any other $E$ given by Kodaira's lemma, and by uniqueness 
we see that $\omega_1$ is smooth off  $E'$, the intersection of the supports of all such $E$.
We claim that if $\alpha=c_1(L)$ for some line bundle $L$, then $E'$ is equal to the null locus of $L$. By Nakamaye's Theorem all we
need to show is that it is equal to the augmented base locus of $L$. If $x\in X$ is a point outside
the augmented base locus, then there exist $H$ an ample divisor and $k, m$ large enough so that
$x$ is not in the base locus of $mL-\frac{m}{k}H$. But this means that $mL-\frac{m}{k}H\sim N$
where $N$ is an effective divisor that doesn't pass through $x$, and moreover the cohomology class of $L-\frac{1}{m}N$ is K\"ahler. So we can take $\ve=\frac{1}{m}$ and $E=N$, and we see that $E'$ is contained in the
null locus of $L$. Conversely, if $x$ belongs to the null locus, then there exists a subvariety
$V$ through $x$ with $\dim V=k$ and $(L^k\cdot V)=0$. Since the potentials for the current $\omega_1$ are continuous, the self-intersection $\omega_1^k$ is a well-defined closed positive current \cite{bt}, which restricts to a nonnegative Borel measure on $V$.
The integral $\int_V \omega_1^k$ is then equal to the cohomological intersection number 
$(L^k\cdot V)$ (see e.g. Corollary 9.3 in \cite{dem2}) which is zero. But if $x$ is not in $E'$ then $\omega_1$ is smooth and K\"ahler
near $x$ and the volume of $V$ with respect to $\omega_1$ would be positive, 
which is a contradiction.

Fix once and for all an $\ve>0$ small enough so that Kodaira's lemma holds. First of all notice that the classes $\alpha_t-\ve E=\kappa_\ve^t$ are all K\"ahler when
$t$ is close to $1$.
Choose a K\"ahler form
$\chi_\ve\in\kappa_\ve$, let $\sigma\in H^0(X,\Op{X}(E))$ be the canonical section, and fix a Hermitian metric $|\cdot|$ on $E$ such
that the following Poicar\'e-Lelong equation holds
\begin{equation}
\omega-\ve[E]=\chi_\ve-\ve\sqrt{-1}\de\db\log|\sigma|,
\end{equation}
where $[E]$ denotes the current of integration on $E$.
Then we have 
$$\beta_t-\ve[E]=\chi_\ve+(\beta_t-\omega)-\ve\sqrt{-1}\de\db\log|\sigma|,$$
and $\chi_\ve^t=\chi_\ve+(\beta_t-\omega)$ is K\"ahler for $t$ close to $1$.
There are smooth functions $\varphi_t$ solutions of
\begin{equation}\label{maeq}
\omega_t^n=(\beta_t+\sqrt{-1}\de\db\varphi_t)^n=\alpha_t^n\Omega,
\end{equation}
where the positive constants $\alpha_t^n$ approach $\alpha^n$ as $t$ goes to $1$,
and $\sup_X \varphi_t=0$. 
We now derive a uniform $L^\infty$ estimate for $\varphi_t$.
Since the Ricci-flat metrics $\omega_t$ have a uniform upper bound on the diameter by Theorem \ref{diameter} and a uniform positive lower bound for the volume $\int_X \omega_t^n$, classical results of Croke \cite{croke}, Li \cite{li} and Li-Yau \cite{LY} give uniform upper bounds for the Sobolev and Poincar\'e constants of $\omega_t$. We temporarily modify the normalization of $\varphi_t$ by requiring that $\int_X \varphi_t\omega_t^n=0$ and we're going to show that
$|\varphi_t|\leq C$. This will then hold for the original $\varphi_t$ as well, with perhaps a bigger constant. We employ a Moser iteration argument in the following way, inspired by \cite{yau1}. For any $p>1$ we compute
\begin{equation}\label{long}\begin{split}
\int_X |\nabla(\varphi_t|\varphi_t|^{\frac{p-2}{2}})|^2_{\omega_t}\omega_t^n&=
\frac{p^2}{4}\int_X |\varphi_t|^{p-2}|\nabla\varphi_t|^2_{\omega_t}\omega_t^n\\
&=
\frac{np^2}{4}\int_X |\varphi_t|^{p-2}\de\varphi_t\wedge\db\varphi_t\wedge\omega_t^{n-1}\\
&\leq \frac{np^2}{4}\int_X |\varphi_t|^{p-2}\de\varphi_t\wedge\db\varphi_t\wedge
\left(\sum_{i=0}^{n-1}\omega_t^{n-1-i}\wedge \beta_t^i\right)\\
&=-\frac{np^2}{4(p-1)}\int_X \varphi_t|\varphi_t|^{p-2}\de\db\varphi_t\wedge
\left(\sum_{i=0}^{n-1}\omega_t^{n-1-i}\wedge \beta_t^i\right)\\
&=\frac{np^2}{4(p-1)}\int_X \varphi_t|\varphi_t|^{p-2}(\beta_t-\omega_t)\wedge
\left(\sum_{i=0}^{n-1}\omega_t^{n-1-i}\wedge \beta_t^i\right)\\
&=\frac{np^2}{4(p-1)}\int_X \varphi_t|\varphi_t|^{p-2}(\beta_t^n-\omega_t^n)\\
&\leq Cp\int_X |\varphi_t|^{p-1}\omega_t^n,
\end{split}
\end{equation}
where we used \eqref{maeq} in the last inequality. 
Using \eqref{long} and the uniform Sobolev inequality for $\omega_t$ and iterating in a standard way (see e.g. \cite{siubook}) we get \begin{equation}\label{moser}
\|\varphi_t\|_{L^\infty}\leq C\left(\int_X|\varphi_t|^2\omega_t^n\right)^{\frac{1}{2}}.
\end{equation}
We then use the uniform Poincar\'e inequality for $\omega_t$ together with \eqref{long} with $p=2$ and with the H\"older inequality and the fact that
the volume $\int_X \omega_t^n$ is bounded above to get
$$\int_X|\varphi_t|^2\omega_t^n\leq C\int_X|\nabla\varphi_t|^2_{\omega_t}\omega_t^n\leq C\int_X|\varphi_t|\omega_t^n\leq C\left(\int_X |\varphi_t|^2\omega_t^n\right)^{\frac{1}{2}},$$
which gives $\int_X|\varphi_t|^2\omega_t^n\leq C$, and so with \eqref{moser} this completes the proof of the $L^\infty$ bound $|\varphi_t|\leq C$.
Notice that such a bound also follows from \cite{eyss}, but our proof is more elementary.

Outside $E$ we have 
$$\beta_t=\chi_\ve^t-\ve\sqrt{-1}\de\db\log|\sigma|,$$
so that the functions $\psi_t=\varphi_t-\ve\log|\sigma|$ solve
\begin{equation}\label{ma}
(\chi_\ve^t+\sqrt{-1}\de\db\psi_t)^n=\alpha_t^n\Omega=e^{F_\ve^t}(\chi_\ve^t)^n
\end{equation}
there, for some appropriate smooth functions $F_\ve^t$, defined on the whole of $X$.
As $t$ approaches $1$, the K\"ahler forms $\chi_\ve^t$ are uniformly bounded in the smooth
topology (with eigenvalues bounded away from $0$ uniformly), and so are the functions $F_\ve^t$.
Yau's second order estimates \cite{yau1} for the Monge-Amp\`ere equation \eqref{ma} give
\begin{equation}\label{max}
\triangle'_t(e^{-A\psi_t}(n+\triangle_t\psi_t))\geq e^{-A\psi_t}\left(-C_1-C_2(n+\triangle_t\psi_t)+(n+
\triangle_t\psi_t)^{\frac{n}{n-1}}\right),
\end{equation}
where $A, C_1$ and $C_2$ are uniform positive constants, $\triangle_t$ is the Laplacian of $\chi_\ve^t$ and $\triangle'_t$
is the Laplacian of $\chi_\ve^t+\sqrt{-1}\de\db\psi_t$. Now notice that on $X\backslash E$ we have
$$e^{-A\psi_t}(n+\triangle_t\psi_t)=|\sigma|^{A\ve}e^{-A\varphi_t}(n+\triangle_t
\varphi_t-\ve\triangle_t\log|\sigma|),$$
and $$\left|\triangle_t\log|\sigma|\right|\leq C,$$
for some uniform constant $C$. Hence the function $e^{-A\psi_t}(n+\triangle_t\psi_t)$ goes to zero
when we approach $E$, and so its maximum will be attained.
The maximum principle applied to \eqref{max} then gives
$$n+\triangle_t\psi_t\leq C e^{A(\psi_t-\inf_{X\backslash E}\psi_t)},$$
on the whole of $X\backslash E$. But noticing that $\inf_{X\backslash E}\psi_t\geq \inf_X \varphi_t -C$
for a uniform constant $C$, and recalling that $|\varphi_t|\leq C_0$, we get
$$n+\triangle_t\varphi_t\leq C+n+\triangle_t\psi_t\leq C(1+|\sigma|^{-A\ve}).$$
This gives uniform interior $C^2$ estimates of $\varphi_t$ and $\psi_t$ on
compact sets of $X\backslash E$. Then the Harnack estimate of Evans-Krylov gives
uniform $C^{2,\gamma}$ estimates, for some $0<\gamma<1$, and a standard bootstrapping 
argument gives uniform $C^{k,\gamma}$ estimates for all $k\geq 2$, on compact sets of $X\backslash E$,
independent of $t<1$. Thus the family $(\varphi_t)$ is precompact $C^{k,\gamma'}(X\backslash E)$ for any $0<\gamma'<\gamma$, and any limit point $\psi$ belongs to
$PSH(X\backslash E,\omega)$, it satisfies 
$$(\omega+\sqrt{-1}\de\db\psi)^n=\alpha^n\Omega$$
on $X\backslash E$, and is bounded near $E$. Hence $\psi$ extends to a bounded
function in $PSH(X,\omega)$ and the above Monge-Amp\`ere equation holds on $X$
because the Borel measure $(\omega+\sqrt{-1}\de\db\psi)^n$ doesn't charge the analytic set $E$. Then by the uniqueness part
of Theorem 2.1 of \cite{eyss},
we must have $\psi=\varphi$. This implies that $\varphi_t\to\varphi$ in
$C^\infty$ on compact sets of $X\backslash E$, and that $\varphi$ is smooth there.
\end{proof}


\section{Examples}\label{examples}
In this section we will give some examples where Theorem \ref{main1} applies. The constructions are well-known and come from algebraic geometry.\\

Let's look at the case $n=2$ first, the case $n=1$ being trivial. The only projective Calabi-Yau surfaces are tori, bi-elliptic, Enriques and
$K3$ surfaces (recall that the Calabi Conjecture has been successfully applied to the study of $K3$ surfaces by Todorov \cite{todorov} and Siu \cite{siu}). If $X$ is a torus and $L$ is a nef and big line bundle on $X$, then $L$ is ample, and so Theorem \ref{main1}
is vacuous in this case. Similarly if $X$ is bi-elliptic, then $X$ is a finite unramified quotient of a torus, so
a nef and big line bundle on $X$ pulls back to a nef and big line bundle on a torus. But this must be ample,
and so the original line bundle is ample too (Corollary 1.2.28 in \cite{laza}) 
and Theorem \ref{main1} is again empty. If $X$ is an Enriques surface,
then $X$ is an unramified $2:1$ quotient of a $K3$ surface, so the study of Ricci-flat metrics on $X$ is reduced
to the case of a $K3$ surface. Finally let's see that there exist projective $K3$s that admit a nef and big line bundle
that is not ample, to which Theorem \ref{main1} applies.
For example let $Y$ be the quotient surface $T/i$ where $T$ is the standard torus $\mathbb{C}^2/\mathbb{Z}^4$
and $i$ is induced by the involution $i(z,w)=(-z,-w)$ of $\mathbb{C}^2$. The surface $Y$ has 16 singular
points, that are rational double points, and is a Calabi-Yau model. Blowing up these 16 points gives a smooth
projective $K3$ surface $X$ (called a Kummer surface), and we can take $L$ to be the pullback of any ample divisor on $Y$. The set $E$, being equal to the null locus of $L$,
is readily seen to be the union of the $16$ exceptional divisors, that are $(-2)$-curves. Then Theorem \ref{main1}
applies, and the limit of smooth Ricci-flat metrics on $X$ with classes approaching $c_1(L)$ is the pullback
of the unique Ricci-flat (actually flat) orbifold K\"ahler metric on $Y$ in the given class. This originally appeared as Theorem 8 in \cite{kt}.
Now we show that conversely all examples of Theorem \ref{main1} on $K3$ surfaces with $\alpha=c_1(L)$ 
are of the form $f:X\to Y$ where $Y$ is an orbifold $K3$ surface, $kL=f^*A$, for some $k\geq 1$ and some $A$ ample divisor on $Y$. Let $X$ be a projective $K3$ surface and $L$ a nef and big line bundle on $X$. 
By Theorem \ref{bpf} we know that
some power $kL$ is globally generated, and we might as well assume that $k=1$. Then the contraction map $f$ of $L$
contracts an irreducible curve $C$ to a point if and only if $C\cdot L=0$.  But since $L\cdot L>0$, the Hodge
Index theorem implies that $C\cdot C<0$.
The long exact sequence in cohomology associated to the sequence
$$0\to \Op{X}(-C)\to\Op{X}\to\Op{C}\to 0,$$
gives that $H^1(X,\mathcal{O}(-C))=0$. Serre duality on the other hand gives $H^2(X,\mathcal{O}(C))=H^0(X,\mathcal{O}(-C))=0$, and $H^1(X,\mathcal{O}(C))=H^1(X,\mathcal{O}(-C))=0$. Riemann-Roch then gives
$$\dim H^0(X,\mathcal{O}(C))=2+\frac{1}{2}C\cdot C,$$ which implies that $C\cdot C$ must be even.
But since $\pi(C)=\frac{C\cdot C}{2}+1$, the virtual genus of $C$, is nonnegative, we see that
$C\cdot C=-2$. This implies that $\pi(C)=0$ and so $C$ is a smooth rational curve with self-intersection $-2$.
Then the point $f(C)$ is a rational double point, and so $Y=f(X)$ is an orbifold $K3$ surface.
Notice that Ricci-flat orbifold metrics on $Y$ exist by \cite{yau1}, \cite{koba}.

Now we turn to examples in dimension $3$. The first one is known as \emph{conifold} in the physics literature \cite{str},
and is described in detail in section 1.2 of \cite{rossi}, for example. Roughly speaking, it is a $3$ dimensional
Calabi-Yau model $Y$ that sits in $\mathbb{P}^4$ as a nodal quintic. It has $16$ singular points, that are nodes and
not of orbifold type. Moreover there exists a small resolution $f:X\to Y$, that is a birational morphism with
$X$ a smooth Calabi-Yau threefold, that is an isomorphism outside the preimages of the nodes, which are $16$
rational curves. If $L$ is the pullback of any ample divisor on $Y$, then $L$ is nef and big on $X$, 
and the limit of smooth Ricci-flat metrics on $X$ with classes approaching $c_1(L)$ is the pullback
of the unique singular Ricci-flat metric on $Y$, which exists by \cite{eyss}. The convergence is smooth on
compact sets outside the union of the $16$ exceptional curves (which is clearly equal to the null locus of $L$). 
There are also other $3$ dimensional examples where the singularities of $Y$
are not isolated: one of these is described in Example 4.6 in \cite{wilson1}, and $Y$ has a curve $C$ of
singularities. Blowing up $C$ gives a Calabi-Yau threefold $X$; if $L$ is the pullback of any ample divisor on $Y$, then the null locus of $L$ is the exceptional divisor $S$ which is a smooth surface ruled over $C$. Again our Theorem \ref{main1} applies, and the convergence is smooth off $S$.

\section{Further directions}\label{direc}
First let us mention an interesting question that arises from Theorem \ref{main1}. We know that
on $X\backslash E$ the Ricci-flat metrics converge smoothly on compact sets to an incomplete Ricci-flat metric $\omega_1$. Its metric completion is a metric space $(X_{\infty},d_{\infty})$. Do the original metrics $(X,\omega_t)$ actually converge to $(X_{\infty},d_{\infty})$ in the Gromov-Hausdorff topology?
We can prove this in the case when $X$ is a $K3$ surface and $\alpha=c_1(L)$. 
In fact, from section \ref{examples} we know that $E$ is a union of $(-2)$-curves and the
contraction map $f:X\to Y$ maps them to orbifold points. The results of \cite{and}, \cite{bkn}, \cite{tian}
give that a subsequence of $(X,\omega_t)$ converges to $Y$ with its orbifold Ricci-flat metric
$d_\infty$ in the Gromov-Hausdorff topology. But on the smooth part of $Y$ we have
that $\omega_1$ and $d_\infty$ coincide, because they are both singular Ricci-flat metrics on the
whole of $Y$. Hence the metric completion of $\omega_1$ is $d_\infty$.

Also, when $X$ admits a birational Calabi-Yau model $Y$, which has a singular Ricci-flat
metric by \cite{eyss}, what is the relation between $Y$ and $X_\infty$? 

There are two possible directions where it would be desirable to extend Theorem \ref{main1}. The first case is when
we look at the whole K\"ahler cone, instead of just the ample cone, and possibly drop the projectiveness assumption.
Suppose $X$ is a compact Calabi-Yau $n$-fold and fix $\omega_0$ a Ricci-flat metric on $X$. 
The N\'eron-Severi space $N^1(X)_\mathbb{R}$ embeds into
$$H^{1,1}_\mathbb{R}(X)=H^2(X,\mathbb{R})\cap H^{1,1}(X),$$
but in general it is a proper subspace (for example a generic projective K3 has $\rho(X)=1<20=\dim H^{1,1}_\mathbb{R}(X)$). 
Inside $H^{1,1}_\mathbb{R}(X)$ we have $\mathcal{K}$, the K\"ahler cone, and its closure $\overline{\mathcal{K}}$, the nef cone. We have that
$$\Kns=\K\cap N^1(X)_\mathbb{R},$$
and similarly for the nef cone.
Given a nonzero class
$\alpha\in\Kb\backslash\K$, and a smooth path $\alpha_t:[0,1]\to\Kb$ that ends at $\alpha$, Yau's Theorem gives a path
$\omega_t$ of Ricci-flat metrics in $\alpha_t$ and we can analyze their behaviour as $t$ approaches $1$. Let's assume  that $\alpha$ is big, which again means that $\alpha^n>0.$
We would like to repeat the construction we did in the algebraic case. There are two main points where we used the assumption that
$X$ was projective and that $\alpha$ belonged to the N\'eron-Severi space: Proposition \ref{semipositive} and Kodaira's lemma. 
We conjecture that the analogue of Proposition \ref{semipositive} still holds, namely we propose the

\begin{conj}
Let $X$ be a compact K\"ahler Calabi-Yau manifold, and $\alpha\in H^{1,1}_\mathbb{R}(X)$
be a class which is nef and big, but not K\"ahler. Then $\alpha$ can be represented by a
smooth $(1,1)$ form $\omega$ which is pointwise nonnegative and which is K\"ahler outside
a proper analytic subvariety $E\subset X$.
\end{conj}
Notice that the proof of this conjecture would have to use the fact the $X$ is Calabi-Yau, since
in general a nef and big class cannot be represented by a smooth nonnegative form (see Example 3.15 in \cite{bb} which is based on \cite{dps2}).
If this conjecture were proved, we could then write $\omega$ as the smooth limit of K\"ahler forms in $\alpha_t$, as in Proposition \ref{semipositive}. 
The correct substitute for Kodaira's lemma would then be given by the theory of closed positive currents: following \cite{paun2}, which relies on the fundamental \cite{dp}, we know that there would exist a modification $\pi:\ti{X}\to X$ such that
$$\pi^*\omega=\ti{\omega}+[E]-\sqrt{-1}\de\db\eta,$$
where $\ti{\omega}$ is a K\"ahler form on $\ti{X}$, $E$ is an effective $\Q$-divisor on $\ti{X}$ and
$\eta$ is quasi-psh, smooth of $E$ and has only log poles along $E$. 
Then we could just work on $\ti{X}$,
and get the same estimates as above, outside $E$, thus proving the K\"ahler analogue of Theorem \ref{main1}.

The second direction is to look at the case when the class $\alpha$ is nef but not necessarily big. Notice that Theorem \ref{diameter} gives a uniform diameter bound in this case. 
A guiding example is
the following: let $X$ be an elliptically fibered $K3$ surface, so $X$ comes equipped with a morphism $f:X\to\mathbb{P}^1$
with fibers elliptic curves. Then the pullback of an ample line bundle on $\mathbb{P}^1$ gives a nef line bundle $L$ on
$X$ with Iitaka dimension $1$. In the case when all the singular fibers of $f$ are of Kodaira type $I_1$, Gross-Wilson have
shown in \cite{gw} that sequences of Ricci-flat metrics on $X$ whose class approaches $c_1(L)$ converge in $C^\infty$
on compact sets of the complement of the singular fibers to the pullback of a K\"ahler metric on $\mathbb{P}^1$.
This metric on $\mathbb{P}^1$ was first studied by McLean \cite{mclean}. In a recent paper, Song-Tian \cite{st} gave a more
direct proof of the result of Gross-Wilson. Moreover they noticed that
McLean's metric satisfies an elliptic equation outside the images of the singular fibers,
namely its Ricci curvature equals the pullback of the Weil-Petersson metric from the moduli space of elliptic curves,
that comes from the variation of the complex structure of the fibers of $f$. 

We believe that in higher dimensions
a similar picture should be true, when $\alpha=c_1(L)$. In this case Conjecture \ref{bpf2} would imply
the existence of a morphism $f:X\to Y$ with connected fibers, where $\dim Y=\kappa(X,L)<n$. Then we expect that
outside a proper subvariety $E\subset X$, a sequence of Ricci-flat metrics with class approaching $\alpha$ should
converge in $C^\infty$ on compact sets of $X\backslash E$ to the pullback of a metric on $Y$. 
It is readily verified that, up to a subsequence, the Ricci-flat metrics converge weakly as currents to the pullback
of a metric on $Y$.
The fibers of $f$ are again Calabi-Yau's, and a computation as in \cite{st} shows that the limit metric on $Y$ will satisfy the same equation as McLean's metric (in this case the potentials of the Ricci-flat metrics have a uniform $C^0$ bound \cite{st}). It might be possible to construct higher-dimensional examples of this behaviour using the results of Section 8 in \cite{fine}, where the equation of McLean's metric appears in his condition (C).

The situation is different when $\alpha$ is not $c_1(L)$, and $X$ possibly not projective.
Then an example of McMullen \cite{ctm} shows that the Ricci-flat metrics can converge smoothly to
zero on an open set of $X$. Also easy examples on tori show that the fibration structure as above
cannot be expected when the limiting class $\alpha$ is not rational. Instead we still expect the Ricci-flat metrics to converge smoothly on compact sets outside a subvariety $E$ to a limit nonnegative form $\omega$,
whose determinant vanishes identically. The kernel of $\omega$ would then define a complex 
foliation with singularities on $X\backslash E$, whose leaves might be dense in $X$. The leaves of the foliation
are always complex submanifolds, but they might not vary holomorphically and the rank of
the foliation might change on different open sets (as in McMullen's example). Notice that if 
the curvature is uniformly bounded, then Ruan's result \cite{ruan} implies that this picture
is basically true and moreover that the 
foliation is holomorphic, so that its rank is constant on a Zariski open set. In McMullen's
example the curvature blows up, and the resulting foliation is not holomorphic, thus showing
that Ruan's result doesn't hold if the curvature is unbounded.

Let us mention that the results of  \cite{bkn}, \cite{bando} also give a description of the behaviour of the Ricci-flat metrics near the singularities, where some bubbling occurs. Unfortunately our methods don't seem to give results of this kind and it would be very interesting to study this in higher dimensions when the limit Calabi-Yau model doesn't have orbifold singularities.

Finally let us notice that some of the results here generalize to the following setting: $X$ is a compact
K\"ahler manifold, and we fix a smooth volume form $\Omega$. If $\alpha_t$ is a path of K\"ahler
classes as in the beginning of this section, then for each $t<1$ Yau's theorem \cite{yau1} gives a unique K\"ahler form
$\omega_t$ in $\alpha_t$ such that $$\omega_t^n=\frac{\alpha_t^n}{\int_X\Omega}\Omega.$$
We can then study the behaviour of the metrics $\omega_t$ as $t$ approaches $1$. If the image of $\alpha_t$ lies in $N^1(X)_\R$ and the limit class
$\alpha$ is nef, big and semiample, then the argument of Theorem \ref{main1} goes through, and we get smooth convergence on compact sets outside a subvariety. Again if $\alpha=c_1(L)$
then this subvariety is the null locus of $L$.


\end{document}